\documentclass[11pt]{article}

\usepackage{amsmath}
\usepackage{amsfonts}
\usepackage{hyperref}
\usepackage{graphicx}
\usepackage{multirow}
\usepackage{xcolor}

\newcommand{\R}{\mathbb{R}}
\newcommand{\N}{\mathbb{N}}
\newcommand{\half}{\frac{1}{2}}
\newcommand{\obs}{\mathrm{obs}}
\newcommand{\pred}{\mathrm{pred}}
\newcommand{\ntime}{n_t}
\newcommand{\nobs}{n_o}
\newcommand{\ndiv}{n_x}
\newcommand{\nred}{n_{\mathrm{red}}}
\newcommand{\ftarget}{f_{\mathrm{target}}}
\newcommand{\xtrial}{x^{\mathrm{trial}}}
\newcommand{\old}{\mathrm{old}}
\newcommand{\accel}{\mathrm{accel}}

\newtheorem{teo}{Theorem}[section]
\newtheorem{lem}{Lemma}[section]

\newcommand{\halmos}{\hfill$\Box$}
\newenvironment{pro}{\noindent\textit{Proof:}}{\halmos}

\DeclareMathOperator*{\Minimize}{Minimize}

\begin{document}

\title{Accelerated derivative-free nonlinear least-squares\\ applied to
  the estimation of Manning coefficients\thanks{This work was
    supported by FAPESP (grants 2013/07375-0, 2016/01860-1, and
    2018/24293-0) and CNPq (grants 302538/2019-4 and 302682/2019-8).}}

\author{
  E. G. Birgin\thanks{Department of Computer Science, Institute of
    Mathematics and Statistics, University of S\~ao Paulo, Rua do
    Mat\~ao, 1010, Cidade Universit\'aria, 05508-090, S\~ao Paulo, SP,
    Brazil. e-mail: egbirgin@ime.usp.br}
  \and
  J. M. Mart\'{\i}nez\thanks{Department of Applied Mathematics,
    Institute of Mathematics, Statistics, and Scientific Computing
    (IMECC), State University of Campinas, 13083-859 Campinas SP,
    Brazil. e-mail: martinez@ime.unicamp.br}}

\date{April 6, 2021}

\maketitle

\begin{abstract}
A general framework for solving nonlinear least squares problems
without the employment of derivatives is proposed in the present paper
together with a new general global convergence theory. With the aim to
cope with the case in which the number of variables is big (for the
standards of derivative-free optimization), two dimension-reduction
procedures are introduced. One of them is based on iterative subspace
minimization and the other one is based on spline interpolation with
variable nodes. Each iteration based on those procedures is followed
by an acceleration step inspired in the Sequential Secant Method. The
practical motivation for this work is the estimation of parameters in
Hydraulic models applied to dam breaking problems. Numerical examples
of the application of the new method to those problems are given.\\

\end{abstract}

\section{Introduction}

Many statistical learning problems require fitting models to large
data sets. Frequently, the number of unknown parameters is not
small. Moreover, for different reasons, derivatives of the functions
that define the model may  not be available, and the sum of squares of
residuals is a natural function to be minimized. These considerations
lead to the problem
 \begin{equation} \label{theprob}
\Minimize \| F(x) ||_2^2 \mbox{ subject to } x \in \Omega \subseteq \R^n,
\end{equation}
where $F:\Omega \to \R^m$.  

Let us define $f(x) = \half \|F(x)\|_2^2$.  For obtaining a quadratic
approximation of~$f(x)$ with the property of being exact if~$f(x)$ is
quadratic, $1 + n + n (n+1) / 2$ evaluations of~$f(x)$ are
needed. However, if the structure $\half \|F(x)\|_2^2$ of $f(x)$ is
used, the same property can be obtained using only~$n+1$ evaluations
of~$F(x)$. Considering that evaluating~$f(x)$ and~$F(x)$ has the same
cost, this seems to be a strong argument to take advantage of the
sum-of-squares structure of~$f(x)$, especially if derivatives are not
available.

Ralston and Jennrich~\cite{ralston} introduced a purely local method
that, at each iteration, minimizes the norm of the linear model that
interpolates~$n+1$ consecutive residuals, providing the first
generalization of the Sequential Secant Method~\cite{barnes,wolfe} to
nonlinear least squares. Zhang, Conn, and Scheinberg~\cite{zcs}
employed different quadratic models for each component of the residual
function in order to define conveniently structured trust-region
subproblems. Therefore, as in \cite{powell2002,powell2013,bobyqa}, at
least~$2n+1$ residual evaluations are computed per iteration. The use
of quadratic models allow these authors to prove not only global
convergence, but also local quadratic convergence under suitable
assumptions~\cite{zhangconn}. The idea of interpolating a different
quadratic for each component of the residual has also been exploited
in the POINDERS software~\cite{wild} with trust-region strategies for
obtaining global convergence. Cartis and Roberts
\cite{cartisroberts2019} introduced a derivative-free Gauss-Newton
method for solving nonlinear least squares problems. At each iteration
of their method, $n+1$ residuals are used to interpolate a linear
model of $F(x)$. The norm of the linear model is approximately
minimized over successive trust regions until sufficient decrease of
the sum of squares is obtained. The points used for interpolation are
updated in order to preserve well-conditioning. With this framework
global convergence and complexity results are proved. All mentioned
methods suffer from a high linear algebra cost per iteration related
to construct the model and find a model's solution; so a natural
idea, explored in the present work, is to apply
dimensionality-reduction techniques. While developing this work, we
became aware of a work of Cartis and Roberts~\cite{cartisroberts2021}
in which this idea is explored. In~\cite{cartisroberts2021}, a method
that performs successive minimizations within random subspaces,
employing the model-based framework based on the Gauss-Newton method
introduced in~\cite{cartisroberts2019}, is introduced.

Having in mind the estimation of parameters in one-dimensional models
that simulate water or mud flow in natural channels, a derivative-free
method for large-scale least-squares problems is introduced in the
present work. Mathematical models for this type of phenomena consist
of partial differential equations with boundary conditions that
simulate flood intensity. The initial conditions for this type of
models are, in general, well known, but the parameters reflecting
density, friction, obstacles, or terrain features must be estimated
from data. We are particularly interested in Manning coefficients.
Manning's coefficients play a crucial role in the correct modeling of
mud or water flow in a natural channel influenced by a flood. In
principle, in perfectly straight channels with constant
cross-sectional area, these coefficients account for velocity
reductions due to friction with the walls or viscosity of the
fluid. In real situations, in which the channel is not straight and
the cross-sectional area is not constant, Manning's coefficients
absorb the information due to these ``irregularities'', which, in
fact, detract the theoretical model from the real situation. The
realistic simulation of a natural channel cannot rely on theoretical
estimates of Manning's coefficients based on physical considerations
linked to ideal situations. Necessarily, such coefficients must be
estimated on the basis of (much or little) available data. This is the
exercise we propose in the present work, for which we use an ideal
situation that allows us to infer the usefulness of the introduced
methods in more realistic situations. Incidentally, data collection in
real cases was dramatically interrupted in 2020 by the outbreak of the
pandemic we are still suffering from. The use of programs whose source
code is not available is frequent in this type of research. For this
reason, we are interested in investigating the behavior of
derivative-free methods to estimate parameters of the models used. The
introduced method combines dimensionality-reduction
techniques~\cite{yaxiang1,yaxiang2} and acceleration steps based on
the sequential secant approach~\cite{barnes,wolfe}.

Acceleration schemes, by means of which, given an iterate $x^k$ and
its predecessors, one obtains a possible (accelerated) better
approximation to the solution may be applied to any of the algorithms
mentioned in the previous paragraph. Let us provide a rough
description of the sequential secant idea applied to nonlinear least
squares problems. Assume that $p \in \{1, 2, \dots n\}$ is given and
$x^0, x^{-1}, \ldots, x^{-p} \in \R^n$ are arbitrary. Given $k = 0, 1,
2,\dots$, we define
\[
s^{k-1} = x^k - x^{k-1}, \dots, s^{k-p} = x^{k-p+1} - x^{k-p},
\]
\[
y^{k-1} = F(x^k) - F(x^{k-1}), \ldots, y^{k-p} = F(x^{k-p+1}) - F(x^{k-p}),
\]
and
\[
S_k = (s^{k-1}, \dots, s^{k-p}) \mbox{ and } Y_k = (y^{k-1}, \dots, y^{k-p}).
\]
The Sequential Secant Method for nonlinear least-squares is defined by
\begin{equation} \label{itessm}
x^{k+1} = x^k - S_k Y_k^\dagger F(x^k),
\end{equation}
where $Y_k^\dagger$ denotes the Moore-Penrose pseudo-inverse of $Y_k$.
Its main drawback is that, according to~(\ref{itessm}), $x^{k+1} -
x^k$ always lies in the subspace generated by $\{ s^{k-1}, \dots,
s^{k-p}\}$. Therefore, all the iterates lie in the affine subspace
that passes through $x^0$ and is spanned by $\{ s^{-1}, \dots, s^{-p}
\}$. This is not a serious inconvenient if $p = n$ and the increments
$s^{k-1}, \dots, s^{k-p}$ remain linearly independent. However, even
when $p = n$, the vectors $s^{k-1}, \dots, s^{k-p}$ may become
linearly dependent and, consequently, all the iterates $x^{k+j}$ would
be condemned to lie in a fixed affine subspace of dimension strictly
smaller than~$n$. For these reasons, the pure Sequential Secant Method
is not appropriate for solving nonlinear least-squares problems
when~$n$ is large and, thus, it is required to maintain $p$ reasonable
small.

Note, however, that, when $m=n$, under suitable assumptions, the
method defined by~(\ref{itessm}) has Q-superlinearly local converge to
a solution of $F(x)=0$, and its R-rate of convergence is the positive
root of $t^{n+1} - t^n - 1 = 0$~\cite{ortegarheinboldt}. When $m=n$,
the problem consists of solving the nonlinear system~$F(x)=0$. This
case has been extensively considered in \cite{birginmartinez2021}. The
drawback pointed out above was overcame in~\cite{birginmartinez2021}
taking auxiliary residual-related directions. The idea of using
residuals as search directions for solving nonlinear systems of
equations have been introduced and exploited
in~\cite{lmr,lacruzraydan,mmps,varadhangilbert} and analyzed from the
point of view of complexity in \cite{grapiglia}. Unfortunately, in
general nonlinear least-squares problems, it is not possible to use
residuals as search directions. Therefore, in the present paper, we
suggest different alternatives for choosing the first trial point
without residual information at each iteration. This is the place
where the dimensionality-reduction techniques place their role --
trial points are computed by minimizing the least-squares function in
a reduced space. Two alternatives are considered. In one of them,
minimizations within small random affine-subspaces are performed. On
the other one, the reduced problem has as variables nodes and values
of a linear spline from which the values of the original variables are
obtained. After the computation of a suitable trial point, we try an
acceleration step using sequential secant ideas.

It is worth mentioning that sequential secant acceleration is closely
connected with Anderson
acceleration~\cite{anderson,brezinski,brezinskiredivo,brs,how,niwalker,rohwedderschneider,wwy}
and quasi-Newton
acceleration~\cite{brezinski,fangsaad,hbdb,mrm}. Moreover, the
sequential secant algorithm is a particular case of a family of secant
methods described in~\cite{ortegarheinboldt} and~\cite{jankowska},
whereas related multipoint secant methods for solving nonlinear
systems and minimization have been introduced
in~\cite{bhd2020,bhd2021,gratton,gmt,scheufelemell,schnabel} and
others.
      
This paper is organized as follows. In Section~\ref{general}, we
introduce a general scheme that applies to derivative-free
optimization (not only nonlinear least-squares) and has the proposed
algorithm for nonlinear least squares as particular case. Global
convergence results for the general scheme are included in this
section. In Section~\ref{algo}, we define the specific algorithm that
we use for derivative-free nonlinear least-problems. In
Section~\ref{numerical}, we present the problem of estimating Manning
coefficients and report numerical experiments. Conclusions and lines
for future research are stated in Section~\ref{conclusions}.\\

\noindent
\textbf{Notation.} The symbol $\|\cdot\|$ will denote an arbitrary
norm.
 
\section{General optimization framework} \label{general}

In this section, we consider the problem
\begin{equation} \label{genprob}
\Minimize f(x) \mbox{ subject to } x \in \Omega,
\end{equation}
where $f: \R^n \to \R$ is arbitrary and $\Omega \subset \R^n$ is
closed and convex.

The following algorithm applies to the solution of
(\ref{genprob}). This algorithm resembles the classical Frank-Wolfe
algorithm \cite{frankwolfe} as, at each iteration, minimizes a linear
function subject to the true constraints of the problem and an
additional constraints that guarantees that the problem is solvable.\\

\noindent
\textbf{Algorithm~\ref{general}.1.} Let $\ftarget \in \R$, $\Delta >
0$, $\gamma \in (0,1)$, a sequence $\{\eta_k\}$ of positive numbers
such that
\begin{equation} \label{sumaetak}
\sum_{k = 0}^\infty \eta_k < \infty,
\end{equation}
and the initial guess $x^0 \in \Omega$ be given. Set $k \leftarrow 0$.

\begin{description}
\item[Step 1.] If $f(x^k) \leq \ftarget$, then terminate the execution
  of the algorithm.
\item[Step 2.] Choose $v^k \in \R^n$ such that $\|v^k\| = 1 $.
\item[Step 3.] Compute $d^k$ as a solution to the suproblem given by
\begin{equation} \label{soldk}
\Minimize \langle v^k, d \rangle \mbox{ subject to } \| d \|
\leq \Delta \mbox{ and } x^k + d \in \Omega.
\end{equation}
\item[Step 4.] Set $\alpha \leftarrow 1$.
\item[Step 5.] Set $\xtrial \leftarrow x^k + \alpha d^k$.
\item[Step 6.] Test the descent condition
\begin{equation} \label{algarmijo}
f(\xtrial) \leq f(x^k) + \eta_k - \gamma \alpha^2 \left[ f(x^k) - \ftarget \right].
\end{equation}
\item[Step 7.] If (\ref{algarmijo}) holds define $\alpha_k = \alpha$,
  compute $x^{k+1} \in \Omega$ such that
\begin{equation} \label{masaun}
   f(x^{k+1}) \leq f( \xtrial),
\end{equation}
 set $k \leftarrow k+1$, and go to Step~1. Otherwise, update $\alpha
 \leftarrow \alpha/2$ and go to Step~5.
\end{description}

\begin{lem} \label{lemafirst}
Assume that $f$ is continuous, $x^k \in \R^n$ is an arbitrary iterate
of Algorithm~\ref{general}.1, and $f(x^k) > \ftarget$. Then, $\xtrial$
and $x^{k+1}$ satisfying~(\ref{algarmijo}) and (\ref{masaun}) are well
defined.
\end{lem}

\begin{pro}
The thesis follows from the continuity of $f$ using that $\eta_k > 0$
and that the successive trials for~$\alpha$ tend to zero.
\end{pro}

\begin{lem} \label{lem31}
Assume that $f$ is continuous and, for all $k\in \N$, we have that
$f(x^k) > \ftarget$. Then,
\begin{equation} \label{produ}
\lim_{k \to \infty} \alpha_k^2 \left[ f(x^k)-\ftarget \right] = 0.
\end{equation}                
Morever, at least one of the following two possibilities takes place:
\begin{equation} \label{lasufi} 
\lim_{k \to \infty} \alpha_k  = 0;
\end{equation}    
or there exists an infinite subset of indices $K_1 \subset \N$ such
that
\begin{equation} \label{limfk1}
\lim_{k \in K_1} f(x^k) = \ftarget. 
\end{equation}
\end{lem}

\begin{pro}
By Lemma~\ref{lemafirst} and the hypothesis, the algorithm generates
an infinite sequence $\{x^k\}$ such that $\{f(x^k)\}$ is bounded
below. Assume that (\ref{produ}) is not true. Then, there exists $c >
0$ such that
\begin{equation} \label{hayc}
\alpha_k^2 [f(x^k) - \ftarget] \geq c
\end{equation}
for infinitely many indices $k \in K_2$. By the convergence of
$\sum_{k = 0}^\infty \eta_k$, there exists $k_1 \in \N$ such that
\[
\eta_k < \gamma c/2
\]  
for all $k \geq k_1$. Then, by (\ref{algarmijo}), (\ref{masaun}), and
(\ref{hayc}), for all $k \in K_2$ such that $k \geq k_1$,
\begin{equation} \label{lagama}
f(x^{k+1}) \leq f(x^k) + \gamma c/2 - \gamma c = f(x^k) - \gamma c/2.
\end{equation}
Let $k_2 \geq k_1$ such that
\[
\sum_{k = k_2}^\infty \eta_k < \gamma c/4.
\]
Then, by (\ref{algarmijo}) and (\ref{masaun}), for all $k \in \N$, $ k
> k_2$ we have that
\begin{equation} \label{cuantosube}
\begin{array}{rcl}
f(x^k) - f(x^{k_2}) &=& [f(x^k) - f(x^{k-1})]+ [f(x^{k-1})-f(x^{k-2})] +
\dots + [f(x^{k_2+1}) - f(x^{k_2})]\\[2mm]
&\leq& \eta_{k-1} + \eta_{k-2} + \ldots + \eta_{k_2} < \gamma c/4.
\end{array}
\end{equation}
Thus, between two consecutive terms (not smaller than $k_2$) of the
sequence $K_2$, by~(\ref{cuantosube}), $f$ increases at most $\gamma
c/4$; but, by~(\ref{lagama}), decreases at least $\gamma c/2$. This
implies that $\lim_{k \to \infty} f(x^k) = -\infty$, which contradicts
the fact that $\{f(x^k)\}$ is bounded below. Therefore, (\ref{produ})
is proved.

Now, if (\ref{lasufi}) does not hold, there exists an infinite set of
indices $K_1$ such that $\alpha_k$ is bounded away from zero. By
(\ref{produ}), we have that (\ref{limfk1}) must take place.
\end{pro}\\

By Lemmas \ref{lemafirst} and~\ref{lem31}, there are three
possibilities for the sequence generated by Algorithm~\ref{general}.1:
\textbf{(i)} The sequence terminates at some~$x^k$ where $f(x^k) \leq
\ftarget$; \textbf{(ii)} The sequence terminates at some~$x^k$ where
$f(x^k) \leq \ftarget + \varepsilon_f$ for a given tolerance
$\varepsilon_f > 0$; and \textbf{(iii)} The sequence $\{\alpha_k\}$
tends to zero.  Possibilities (i) and (ii) are symptoms of success of
the algorithm. Possibility (iii) cannot be discarded since, given
$\varepsilon_f>0$, $f(x)$ may be bigger than $\ftarget +
\varepsilon_f$ for every $x \in \Omega$. Therefore, the implications
of $\alpha_k \to 0$ need to be analyzed. With this purpose, in the
following lemma we need to assum differentiability of the
function~$f$.

\begin{lem} \label{lema82}
Assume that $f$ admits continuous derivatives for all $x$ in an open
set that contains~$\Omega$ and $\{x^k\}$ is generated by
Algorithm~\ref{general}.1. Assume that $\{f(x^k)-\ftarget \}$ is
bounded away from zero, $x_* \in \R^n$, and $\lim_{k \in K_1} x^k =
x_*$. Then, the sequence $\{d^k\}_{k \in K_1}$ admits at least one
limit point and, for every limit point $d$ of $\{d^k\}_{k \in K_1}$, we
have that
\begin{equation} \label{ortogo1}
\langle \nabla f(x_*), d \rangle \geq 0.
\end{equation}
\end{lem}            

\begin{pro}
By Lemma~\ref{lem31},
\begin{equation} \label{alfaacero}
\lim_{k \to \infty} \alpha_k = 0.
\end{equation}
Since the first trial value for $\alpha_k$ at each iteration is $1$,
(\ref{alfaacero}) implies that 
\[
\lim_{k \to \infty} \alpha_{k,+}  = 0, 
\]
and, for all $k$ large enough, 
\[
f(x^k + \alpha_{k,+} d^k) > f(x^k) + \eta_k - \gamma \alpha_{k,+}^2
\left[ f(x^k)-\ftarget \right].
\]
So, since $\eta_k > 0$, 
\[
\frac{f(x^k + \alpha_{k,+} d^k ) - f(x^k)}{\alpha_{k,+}} > - \gamma
\alpha_{k,+} \left[ f(x^k)-\ftarget \right]
\]
for all $k$ large enough.  Thus, by the Mean Value Theorem, there exists
$\xi_{k,+} \in [0, \alpha_{k,+}]$ such that
\begin{equation} \label{notarmijo7}
\langle \nabla f(x^k + \xi_{k,+} d^k), d^k \rangle > -\gamma
\alpha_{k,+} \left[ f(x^k) - \ftarget \right]
\end{equation}
for all $k$ large enough. 

Since $\|d^k\| \leq \Delta$ for all $k$, we have that $\{d^k\}_{k\in
  K_1}$ admits at least one limit point. Let $d$ be an arbitrary limit
point of $\{d^k\}_{k\in K_1}$ and let $K_2 \subseteq K_1$ such that
\begin{equation} \label{lik3}
\lim_{k \in K_2} d^k = d
\end{equation}
and $\|d\| \leq \Delta$.  
By continuity, since $\lim_{k \in K_2} x^k = x_*$ we have that
\begin{equation} \label{flim}
\lim_{k \in K_2} f(x^k) = f(x_*).
\end{equation}
Then, taking limits for $k \in K_2$ in both sides
of~(\ref{notarmijo7}), by~(\ref{alfaacero}), (\ref{lik3}),
(\ref{flim}), and the fact that (\ref{alfaacero}) implies $\lim_{k \in
  K_2} \xi_{k,+} = 0$, we get
\[
\langle \nabla f(x_*), d \rangle \geq 0
\]
as we wanted to prove. 
\end{pro}

\begin{lem} \label{lema83}
Assume that $f$ admits continuous derivatives for all $x$ in an open
set that contains~$\Omega$ and $\{x^k\}$ is generated by
Algorithm~\ref{general}.1. Assume that $\{f(x^k)-\ftarget \}$ is
bounded away from zero, $x_* \in \R^n$, and there exists $K_1 \subset
\N$ such that $\lim_{k \in K_1} x^k = x_*$ and
\begin{equation} \label{unitaria}
\lim_{k \in K_1} \left\| v^k - \frac{\nabla f(x^k)}{\|\nabla f(x^k)\|}
\right\| = 0.
\end{equation}     
Then, for all $d \in \R^n$ such that $\|d\|\leq \Delta$ and $x_* + d
\in \Omega$ we have that
\[
\langle \nabla f(x_*), d \rangle \geq 0.
\]
\end{lem}            

\begin{pro}
If $\nabla f(x_*)=0$, we are done; so we assume $\nabla f(x_*) \neq 0$
from now on. By Lemma~\ref{lema82}, there exists $K_2 \subseteq K_1$
and $\bar d \in \R^n$ such that
\begin{equation} \label{dbarlim}
\lim_{k \in K_2} d^k = \bar d
\end{equation}
and
\begin{equation} \label{ortogodos}
\langle \nabla f(x_*), \bar d \rangle \geq 0;
\end{equation}
and, since $\nabla f(x_*) \neq 0$, \eqref{ortogodos} implies
\begin{equation} \label{visigodos}
\left\langle \frac{\nabla f(x_*)}{\|\nabla f(x_*)\|}, \bar d \right\rangle \geq 0.
\end{equation}  
Given $\varepsilon>0$, by~\eqref{unitaria}, \eqref{dbarlim}, and the
continuity of $\nabla f$, \eqref{visigodos} implies that
\begin{equation} \label{godos}
\langle   v^k, d^k \rangle \geq - \varepsilon
\end{equation}
for all $k \in K_2$ large enough. Since, by the definition of
Algorithm~\ref{general}.1, $d^k$ is a solution to~\eqref{soldk},
\eqref{godos} implies that
\begin{equation} \label{vdmayo}
\langle   v^k, d \rangle \geq - \varepsilon 
\end{equation}    
for all $d \in \R^n$ such that $x^k + d \in \Omega$ and $\|d\| \leq
\Delta $ and all $k \in K_2$ large enough.

Consider the problem
\begin{equation} \label{sublim}
\mbox{ Minimize } \left\langle \frac{\nabla f(x_*)}{\|\nabla f(x_*)\|}
, d \right\rangle \mbox{ subject to } \| d \| \leq \Delta \mbox{ and }
x_* + d \in \Omega
\end{equation}
that, by compacity, admits a solution $d_*$; and suppose, by
contradiction, that
\begin{equation} \label{labelest}
\left\langle \frac{\nabla f(x_*)}{\|\nabla f(x_*)\|} , d_*
\right\rangle = - c < 0.
\end{equation}
Therefore,
\[
\left\langle \frac{\nabla f(x_*)}{\|\nabla f(x_*)\|} , x_* + d_* - x_*
\right\rangle = - c < 0.
\]  
This implies, by (\ref{unitaria}), that
\begin{equation} \label{labost}
\left\langle v^k , x_* + d_* - x^k \right\rangle \leq - c/2 < 0
\end{equation}
for $k \in K_2$ large enough. Let us write $\tilde d^k = x_* + d_* -
x^k$. Since $x_* + d_* \in \Omega$, we have that $x^k + \tilde d^k \in
\Omega$. If, for some $k \in K_2$ large enough, we have that $\|
\tilde d^k \| \leq \Delta$, taking $\varepsilon < -c/2$, we get a
contradiction between~(\ref{labost}) and~(\ref{vdmayo}). This
contradiction comes from the assumption~(\ref{labelest}), which, as a
consequence, is false, completing the proof.

We now consider the case in which $\| \tilde d^k \| > \Delta$ for all
$k \in K_2$ large enough. Since
\[
\|\tilde{d}^k \| = \|x_* + d_* - x^k\| \leq \|d_*\| + \|x_* - x^k\|
\leq \Delta + \|x_* - x^k\|,
\]
defining
\begin{equation} \label{dchapeu}
\hat d^k =  \frac{ \tilde d^k}{1 + \|x^k - x_*\|/\Delta},
\end{equation}
we have that $\| \hat d^k \| \leq \Delta$ and, by the convexity of
$\Omega$, $x^k + \hat d^k \in \Omega$. Moreover, by
(\ref{labost}), since $\lim_{k \in K_2} x^k = x_*$, 
\begin{equation} \label{ce4}
\left\langle v^k, \hat d^k \right\rangle \leq - c/4 < 0
\end{equation}
for $k \in K_2$ large enough. Taking $\varepsilon < -c/4$, we get the
contradiction between (\ref{ce4}) and (\ref{vdmayo}) and the proof is
complete.
\end{pro}

\begin{teo} \label{teo84}
Assume that $f$ admits continuous derivatives for all $x$ in an open
set that contains~$\Omega$ and $\{x^k\}$ is generated by
Algorithm~\ref{general}.1. Assume that the level set defined by
$f(x^0) + \eta$ is bounded, where $\eta = \sum_{k=0}^\infty \eta_k$,
and there exists an infinite sequence of indices $K_1$ such that
(\ref{unitaria}) holds. Then, given $\varepsilon > 0$, either exists
an iterate~$x^k$ such that $f(x^k) \leq \ftarget + \varepsilon$ or
there exists a limit point $x_*$ of $\{x^k\}$ such that $\langle
\nabla f(x_*), d \rangle \geq 0$ for all $d$ such that $x_* + d \in
\Omega$.
\end{teo}

\begin{pro}
Since the level set defined by $f(x^0)+\eta$ is bounded, the sequence
$\{x^k\}_{k \in K_1}$ admits a limit point. Therefore, the thesis
follows from Lemma~\ref{lema83}.
\end{pro} \\

\noindent
\textbf{Remark.} Let us show that Assumption~(\ref{unitaria}) is
plausible. With this purpose, assume that it does not hold. Then,
there exists $\varepsilon > 0$ such that for all $k \in \N$,
\[
\left\| v^k -  \frac{\nabla f(x^k)}{\|\nabla f(x^k)\|} \right\| > \varepsilon.
\]
Clearly, if we choose randomly the vectors $v^k$ in the unitary
sphere, the probability of this event is zero. Therefore,
assumption~(\ref{unitaria}) holds with probability 1.

\section{Practical algorithm for nonlinear least-squares} \label{algo}

In this section, we are interested in the application of
Algorithm~\ref{general}.1 to large scale nonlinear least squares
problems of the form
\begin{equation} \label{probnls}
\Minimize \half \| F(x) \|_2^2,
\end{equation}
where $F:\R^n \to \R^m$ and $\| \cdot \|_2$ is the Euclidean
norm. Consequently, we define
\begin{equation} \label{funls}
f(x) = \half \|F(x)\|_2^2.
\end{equation}

The proposed algorithm for solving~(\ref{probnls}) is a particular
case of Algorithm~\ref{general}.1 for the case $\Omega=\R^n$, but it
includes two additional features: minimizations in reduced spaces and
acceleration steps.\\

\noindent
\textbf{Algorithm~\ref{algo}.1.} Let $\ftarget \in \R$, $\Delta > 0$,
$\gamma \in (0,1)$, a sequence $\{\eta_k\}$ of positive numbers such
that
\[
\sum_{k = 0}^\infty \eta_k < \infty,
\]
and the initial guess $x^0 \in \R^n$ be given. Set $k \leftarrow 0$.

\begin{description}
\item[Step 1.] If $f(x^k) \leq \ftarget$, then terminate the execution
  of the algorithm.
\item[Step 2.] Compute $\xtrial \in \R^n$ by means of a Reduction
  Algorithm.
\item[Step 3.] Test the descent condition
  \begin{equation} \label{armijored}
    f(\xtrial) \leq f(x^k) + \eta_k - \gamma \left[ f(x^k) - \ftarget \right].
  \end{equation}
  If $ \xtrial \neq x^k$ and~(\ref{armijored}) holds, set $d^k =
  \xtrial - x^k$, $\alpha_k = 1$, $v^k = d^k / \|d^k\|$, and go to
  Step~9.
\item[Step 4.] Choose $v^k \in \R^n$ such that $\|v^k\| = 1$.
\item[Step 5.] Compute $d^k$ as a solution to the subproblem given by
  \[
  \Minimize \langle v^k, d \rangle \mbox{ subject to } \| d \| \leq \Delta.
  \]
\item[Step 6.] Set $\alpha \leftarrow 1$.
\item[Step 7.] Set $\xtrial \leftarrow x^k + \alpha d^k$.
\item[Step 8.] Test the descent condition
  \begin{equation} \label{armijoblue}
    f(\xtrial) \leq f(x^k) + \eta_k - \gamma \alpha^2 \left[ f(x^k) - \ftarget \right].
  \end{equation}
  If (\ref{armijoblue}) does not hold, update $\alpha \leftarrow
  \alpha/2$ and go to Step~7. Otherwise, set $\alpha_k = \alpha$.
\item[Step 9.] Compute, by means of the Acceleration Algorithm,
  $x^{k+1} \in \R^n$ such that
  \[
  f(x^{k+1}) \leq f(\xtrial).
  \]
  Set $k \leftarrow k+1$ and go to Step~1.
\end{description}

Lemmas \ref{lemafirst}, \ref{lem31}, \ref{lema82}, \ref{lema83}, and
Theorem~\ref{teo84} hold for Algorithm~\ref{algo}.1 exactly in the
same way as they do for Algorithm~\ref{general}.1. In order to
complete the definition of Algorithm~\ref{algo}.1, we now introduce
two possible Reduction Algorithm and an Acceleration Algorithm.  Both
Reduction algorithms employ BOBYQA~\cite{bobyqa} for minimizing $f(x)$
over manifolds of moderate dimension.

\subsection{Affine-subspaces-based Reduction Algorithm} \label{affine}

In this Reduction Algorithm the manifold over which we minimize $f(x)$
at each iteration is an affine subspace.  At iteration $k$, we
consider an affine transformation ${\cal T}_k: \R^{\nred} \to \R^n$,
with $\nred \leq n$, given by ${\cal T}_k(d) := x^k + M_k d$, where
$M_k \in \R^{n \times \nred}$ is a matrix with random (with uniform
distribution) elements $m_{ij} \in [-1,1]$. So the problem to be
solved at iteration $k$ is given by
\begin{equation} \label{subspaces}
\Minimize_{d \in \R^{\nred}} \| F({\cal T}_k(d)) \|_2^2.
\end{equation}
The natural initial guess for this problem is given by~$d=0$.  The
minimization on small-dimensional subspaces has been used
in~\cite{yaxiang1,yaxiang2}. Moreover, in the context of
derivative-free optimization, it has been recently employed in
\cite{cartisroberts2021}.

\subsection{Linear-interpolation-based Reduction Algorithm} \label{splin}

At iteration $k$, we consider a linear-spline-based transformation
${\cal S}_k: \R^{\nred} \to \R^n$, with $\nred \leq n$, where $\nred =
2 \kappa + 2$ for some $\kappa \ge 0$. Variables of the reduced model
are $p_1, \dots, p_{\kappa}$, and $v_0, v_1, \dots, v_{\kappa},
v_{\kappa+1}$, with $0 \leq p_j \leq 1$ for $j=1,\dots,\kappa$. Define
$p_0=0$ and $p_{\kappa+1}=1$. If $p_{j_1}=p_{j_2}=\dots$ for $j_1 \neq
j_2 \neq \dots$, then redefine $v_{j_1}, v_{j_2}, \dots$ as their
average. Let $\bar p_0 < \dots < \bar p_{\bar \kappa+1}$ (with $\bar
\kappa \leq \kappa$) be a permutation of $p_0, \dots, p_{\kappa+1}$ in
which repeated values were eliminated and let $\bar v_0, \dots, \bar
v_{\bar \kappa +1}$ be the corresponding (reordered) values. We define
a piecewise linear function $L:[0,1] \to \R$ such that $L(\bar
p_j)=\bar v_j$ for $j=0,\dots,\bar \kappa+1$. The transformation
${\cal S}_k$ is given by ${\cal
  S}_k(v_0,\dots,v_{\kappa+1},p_1,\dots,p_{\kappa}) := x^k + d(v,p)$,
where $[d(v,p)]_i = L((i-1)/(n-1))$ for $i=1,\dots,n$. So the
\textit{bound constrainted} problem to be solved at iteration $k$ is
given by
\begin{equation} \label{splines}
\Minimize_{(v,p) \in \R^{n_{\mathrm{red}}}} \| F({\cal S}_k(v,p)) \|_2^2
\mbox{ subject to } 0 \leq p_j \leq 1 \mbox{ for } j=1,\dots,\kappa.
\end{equation}
As initial guess, we consider $v=0$ and $p$ with random (with uniform
distribution) components $p_j \in [0,1]$ for $j=1,\dots,\kappa$.

\subsection{Acceleration Algorithm} \label{secaccel}

We adopt a Sequential Secant approach for defining the
acceleration. The scheme, that generalizes the one adopted in
\cite{birginmartinez2021} for solving nonlinear sytems of equations,
is as follows.

\begin{enumerate}

\item If $k = 0$, then define $x^{k+1} = \xtrial$.

\item If $k > 0$, then choose $k_{\old} \in \{0, 1, \dots, k-1\}$,
\[
s^j = x^{j+1} - x^j \mbox{ for all } j = k_{\old}, \dots, k-1,
\]
\[
s^k = \xtrial - x^k, 
\]
\[
y^j = F(x^j + s^j) - F(x^j) \mbox{ for all } j = k_{\old}, \dots, k.  
\]
\[
S_k = (s^{k_{\old}}, \dots, s^k),
\]
\[
Y_k = (y^{k_{\old}}, \dots, y^k),
\]
\[
x^k_{\accel} = x^k - S_k Y_k^\dagger F(x^k).
\]

\item If $f(x^k_{\accel}) \leq f(\xtrial)$, then define $x^{k+1} =
  x^k_{\accel}$. Otherwise, define $x^{k+1} = \xtrial$.

\end{enumerate}

This algorithm differs from the plain acceleration scheme defined
in~\eqref{itessm} in a very substantial way. In~\eqref{itessm}, the
definition of $x^{k+1}$ depends only on the previous iterates and lies
in the affine subspace determined by them. Therefore,
in~\eqref{itessm}, the successive iterates do not escape from a fixed
$p$-dimensional affine subspace, where $p$ is the number of previous
iterates that contribute to the acceleration process. On the contrary,
here, we define $x^{k+1}$ as the possible result of an acceleration
that includes the trial point $\xtrial$ which, in principle, does not
belong to any pre-determined affine subspace. As a consequence, the
accelerated point has the chance of exploring the whole domain in a
more efficient way.

\section{Estimation of Manning coefficients in the Saint-Venant equation} \label{numerical}

In the present work, we are interested in the estimation of parameters
in one-dimensional models that simulate water or mud flow in natural
channels. The presence of extreme boundary conditions can be a
consequence of upstream levee breakage, a subject that is studied in
the context of the interdisciplinary research anf action group CRIAB
(acronym for ``Dams Conflicts, Risks and Impacts'' in Portuguese) at
the University of Campinas.  The initial conditions for this type of
models are, in general, well known, but the parameters reflecting
density, friction, obstacles, or terrain features must be estimated
from data.  Mathematical models for this type of phenomena consist of
partial differential equations with boundary conditions that simulate
flood intensity. The use of programs whose source code is not
available is frequent in this type of research. For this reason we are
interested in investigating the behavior of derivative-free methods to
estimate parameters of the models used.

For simplicity, in this study we assume that the phenomenon we are
interested in is well represented by the Saint-Venant equations
\cite{saintvenant}. More sophisticated tools are beyond the scope of
the present work. The Saint-Venant equations
\begin{equation} \label{consmassa}
A_t + Q_x = 0
\end{equation}          
and
\begin{equation} \label{quamov}
Q_t + (Q \, V)_x + g \, A \, \widehat z + \frac{\xi \, P \, V \, |V|}{8} = 0.
\end{equation}
simulate the evolution of mean velocity, wetted cross-sectional area,
depth, and flow in a one-dimensional channel. In \eqref{consmassa}
and~\eqref{quamov}, $A = A(x,t)$ is the wetted cross sectional area at
position $x$ and time $t$; $V = V(x,t)$ is the mean velocity; $Q =
Q(x,t) = A(x,t) \, V(x,t)$ is the flow rate; $P = P(x,t)$ is the
wetted perimeter, that is, the perimeter enclosing the wetted area
taking away the air contact surface; $g$ is the acceleration of
gravity, approximately $9.8 m/s^2$; $\widehat z = \widehat z(x,t) =
z_x / (1+(z_x)^2)$, where $z=z(x,t)=h(x,t) + z_b(x)$, $h(x,t)$ is the
maximum channel depth at point $x$ and time $t$, and $z_b(x)$ is the
vertical coordinate of the channel bottom at point $x$ (therefore,
$z_x=h_x+(z_b)_x$); and $\xi = \xi(x)$ is the adimensional Manning
coefficient whose estimation using data is the subject of the present
study. The estimation of Manning coefficients is a very hard problem
related with the simulation of floods in natural
channels~\cite{djy}. In the present work, we adopt that (a) Manning
coefficients vary at different points of the channel but are invariant
in time and (b) the best estimation of Manning coefficients is the one
that provides the best predictions of streams in a period of time.

We assume that the channel under consideration extends
one-dimensionally from $x=x_{\min}$ to $x=x_{\max}$. The boundary
condition on the left ($x_{\min}$) simulates a flow rate that grows
linearly from~$8.245$ $m^3/s$ to~200 $m^3/s$ in $1{,}200$ seconds a
decreases to the initial flow rate between $1{,}200$ seconds and
$3{,}600$ seconds, remaining stationary thereafter. The initial depth
is $1.2$ meters.  The second derivatives of other state variables are
assumed to be zero both in $x=x_{\min}$ and in $x =x_{\max}$. We
consider that wetted cross-sectional areas and velocities are measured
between times $t=t_{\min}$ and $t=t_{\max}^{\obs}$, at equally spaced
points in the interval $[x_{\min}, x_{\max}]$.  The physical
characteristics of this channel were taken from~\cite{porto}
and~\cite{grafaltinakar}.
Synthetic data were created with $x_{\min}=0$ meters, $\Delta x = 6$
meters, $x_{max} \in \{ 3{,}000, 3{,}600, \dots, 9{,}000 \}$ meters,
$t_{\min}=0$ seconds, and $\Delta t = 0.1$ seconds. The value of
$t_{\max}^{\obs}$, maximum observation time, was subject to
experimentation; while $t_{\max}^{\pred}$, maximum time for
prediction, was set to $t_{\max}^{\pred}=3{,}600$ seconds. The
transversal area was considered to be rectangular with a width of 5
meters. We assumed that the true value for the adimensional Manning
coefficients at the discretizated space points is $0.0366$ plus a
random uniform perturbation of up to~1\%. We set $(z_b)_x =
0.001$. For the purposes of this research, we found it satisfactory to
solve the Saint-Venant equations by finite differences using a
Lax-Friedrichs type scheme~\cite{leveque} with artificial diffusion
coefficient equal to~0.9.

The considerations above lead to a problem of the form
\begin{equation} \label{lsprob}
\Minimize_{\xi \in \R^{\ndiv}} \sum_{i=1}^{\ntime} \sum_{j=0}^{\ndiv} \sum_{k=1}^2
\left( y(\xi,t_i,x_j,k) - y_{ijk}^{\obs} \right)^2,
\end{equation}
where $x_j = x_{\min} + j \Delta x$ for $j=0,\dots,\ndiv$, $t_i =
t_{\min} + i \Delta t$ for $i=1,\dots,\ntime$ and $\ndiv$, $x_{\min}$,
$\Delta x$, $\ntime$, $t_{\min}$, and $\Delta t$ are given. When
$k=1$, $y_{ijk}^{\obs}$ ($i=1,\dots,\ntime$, $j=0,\dots,\ndiv$)
corresponds to a given observation of transversal area; while, when
$k=2$, it corresponds to a given observed velocity. The problem
has~$\ndiv$ unknowns and $2 \ntime (\ndiv+1)$ terms in the
summation. $y(\xi,t_{\min},x_j,k)$ does not depend on~$\xi$ and it
assumed to be known for $k=1,2$ and $j=0,\dots,\ndiv$; while the given
values of $\Delta x$ and $\Delta t$ are such that, if $\xi$ and
$y(\xi,t_i,x_j,k)$ for $j=0,\dots,\ndiv$ are known, then
$y(\xi,t_{i+1},x_j,k)$, for $j=0,\dots,\ndiv$, may be computed in
finite time. In a generalization to~(\ref{lsprob}), it is assumed that
most of the observations are not available and, then, (\ref{lsprob})
is substituted with
\begin{equation} \label{lslost}
\Minimize_{\xi \in \R^{\ndiv}} \sum_{\{ (i,j,k) \in S \}}
\left( y(\xi,t_i,x_j,k) - y_{ijk}^{\obs} \right)^2,
\end{equation}
where $S \subseteq \widehat S$ is given, $\widehat S = \{ (i,j,k)
\;|\; i=1,\dots,\ntime, \; j=0,\dots,\ndiv, \; k=1,2 \}$, and $|S| \ll
|\widehat S|$. For further reference, we denote by $\nobs = |S|$ the
number of available observations. Note that, if $S=\widehat S$, then
$\nobs = 2 \ntime (\ndiv+1)$; while if, for example, only 10\% of the
observations are available, then we have $\nobs = 0.2 \ntime
(\ndiv+1)$.

Approximately solving~(\ref{lslost}) provides a value~$\bar \xi$; and
this value~$\bar \xi$ is then used to predict that
\[
y_{ijk}^{\obs} \approx y(\bar \xi,t_i,x_j,k) \mbox{ for } (i,j,k) \in \widehat S^+,
\]
where $t_{\max}^{\obs} = t_{\min} + n_t \Delta t$ is the largest time
instant at which observations considered in~(\ref{lslost}) were
collected, $t_{\max}^{\pred} > t_{\max}^{\obs}$, and $\widehat S^+ =
\{ (i,j,k) \;|\; t_{\max}^{\obs} < t_i \leq t_{\max}^{\pred}, \;
j=0,\dots,\ndiv, \; k=1,2 \}$ represents the set of indices of the
$y_{ijk}^{\obs}$, \textit{not yet observed}, whose predicted value is
given by $y(\bar \xi,t_i,x_j,k)$.

\subsection{Numerical results}

We implemented Algorithm~\ref{algo}.1, together with the two Reduction
Algorithms (Sections~\ref{affine} and~\ref{splin}) and the
Acceleration Algorithm (Section~\ref{secaccel}) in Fortran~90. All
tests were conducted on a computer with a 3.4 GHz Intel Core i5
processor and 8GB 1600 MHz DDR3 RAM memory, running macOS Mojave
(version 10.14.6). Code was compiled by the GFortran compiler of GCC
(version 8.2.0) with the -O3 optimization directive enabled. In the
rest of this section, mainly in figures and tables,
Algorithm~\ref{algo}.1 is sometimes referred to as SESEM, that stands
for ``Sequential Secant Method''. Based on~\cite{birginmartinez2021}
and on preliminar numerical experiments, we set $\gamma=10^{-4}$,
$\eta_k=$ for $k=0,1,\dots$, and $\Delta=10$ in
Algorithm~\ref{algo}.1, and $p=1000$, i.e.\ $k_{\old} = \max\{ 0, k -
p \}$, in the Acceleration Algorithm.

In section~\ref{nr1}, we aim to determine (i) the amount of
observations (starting at $t_{\min}=0$ and at intervals $\Delta t =
0.1$ seconds), determined by the maximum observation time
$t_{\max}^{\obs}$, and (b) the precision of the optimization process
that are required to recover Manning coefficients~$\xi$ suitable for
making predictions up to $t_{\max}^{\pred}=3{,}600$
seconds. Sections~\ref{nr2} and~\ref{nr3} are related to the
calibration and analysis of the proposed method. In Section~\ref{nr2},
the dimension of the subproblem solved at each iteration is
determined; while in Section~\ref{nr3} the influence of the
Acceleration Algorithm in the overall process is observed. In
Section~\ref{nr4}, a set of instances of increasing size, mimic the
the size of real-life instances, is solved. Section~\ref{nr5} presents
the behavior of the solvers BOBYQA~\cite{bobyqa} and DFBOLS~\cite{zcs}
in the set of considered instances.

\subsubsection{Choice of a tolerance that leads to acceptable solutions} \label{nr1}

Given data coming from observations, we seek to estimate the Manning
coefficients by means of which the Saint-Venant equations produce the
best reproduction of data.  In real cases, we are tempted to believe
that an accuracy of around 10\% in the prediction of depths and
velocity is sufficiently good and that more accurate reproduction is
not justified since observation and modeling errors may be, many
times, of that order.  However, we have no guarantees about the
quality of predictions for data that are not available yet; and it can
be argued that, although an excessive precision in the available data
has no effect in the reproduction of these data, it may have a
significant effect in the reproduction of observations that are not
available yet.  Therefore, it is sensible to test our inversion
procedure not only up to the precision compatible with observation and
modeling errors but also with moderate higher precisions.
 
Assume that an iterative optimization process is applied
to~(\ref{lslost}) to compute~$\bar \xi$; and that this process stops
when it finds~$\bar \xi$ satisfying
\begin{equation} \label{stopcrit}
\left[ \sum_{\{ (i,j,k) \in S \}}
\left( y(\bar \xi,t_i,x_j,k) - y_{ijk}^{\obs} \right)^2 \right] \leq \epsilon
\left[ \sum_{\{ (i,j,k) \in S \}} \left( y_{ijk}^{\obs} \right)^2 \right],
\end{equation}
where $\epsilon>0$ is a given tolerance. Of course, $\bar \xi$ depends
on~$\epsilon$ and on the problem data. In particular, $\bar \xi$ depends
on the set of available observations $S$, that depends
on~$t_{\max}^{\obs}$. Assume that, after computing~$\bar \xi$,
$t_{\max}^{\pred} > t_{\max}^{\obs}$ is chosen and observations
$y_{ijk}^{\obs}$ with $(i,j,k) \in S^+ \subseteq \widehat S^+$ become
available. We define that, for the given $S^+$ and $t_{\max}^{\pred}$,
$\bar \xi$ is \textit{acceptable} if we have that
\begin{equation} \label{accep}
\eta(\bar \xi) := \frac{\sum_{\{ (i,j,k) \in S \cup S^+ \}} \left( y(\bar
  c,t_i,x_j,k) - y_{ijk}^{\obs} \right)^2}{\sum_{\{ (i,j,k) \in S \cup
    S^+ \}} \left( y_{ijk}^{\obs} \right)^2} \leq 10^{-4}.
\end{equation}

\begin{figure}[ht!]
\begin{center}
\input{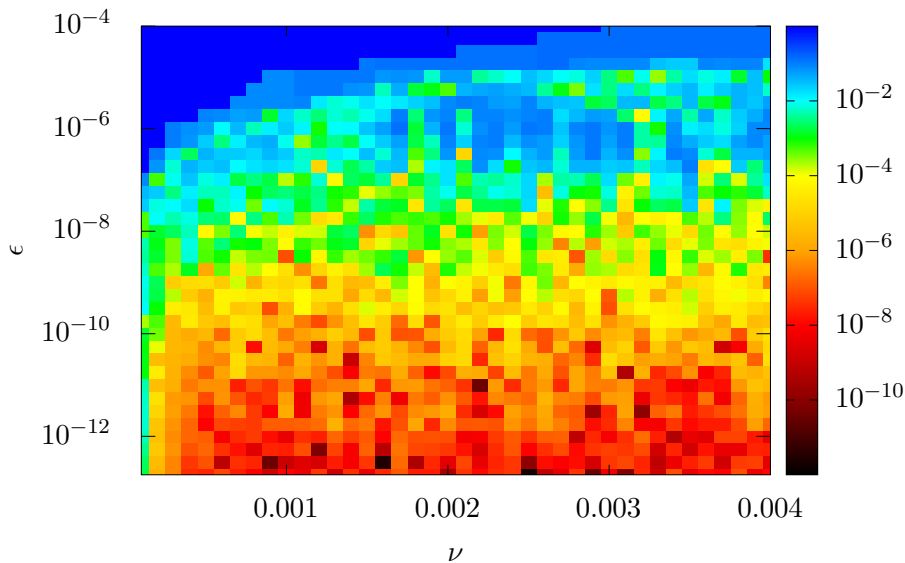}
\end{center}
\caption{Acceptability of solutions $\bar \xi(\nu,\epsilon)$ to
  instances with varying number of observations $\ntime(\nu)$ solved
  with varying tolerances~$\epsilon$. Hot colors show that acceptable
  solutions can be computed when the number of observations is larger
  than the number of unknowns and the tolerance to stop the
  optimization process is tight (smaller than $10^{-9}$).}
\label{fig1}
\end{figure}

Let the problem data $\ndiv$, $x_{\min}$, $\Delta x$, $t_{\min}$, and
$\Delta t$ (note that $\ntime$ is missing here) be given and assume
that an instant $t_{\max}^{\pred}$ is chosen. The question is: Which
are the number of observations~$\nobs$ and the optimization
tolerance~$\epsilon$ that make the computed $\bar \xi$ to be
acceptable?  We aim to answer this question empirically considering a
typical instance of~(\ref{lslost}) with $\ndiv=500$, $x_{\min}=0$,
$\Delta x=6$, $t_{\min}=0$, $\Delta t = 0.1$, and $S$ randomly chossen
in such a way that $\nobs = |S| \approx 0.1 (2 \ntime (\ndiv+1))$,
i.e.\ assuming that approximately 90\% of the observations are not
available. (Units of measure are meters for space and seconds for
time.) Setting $t_{\max}^{\pred} = 3{,}600$ and varying a constant
$\nu \in \{ 10^{-4}, 2 \times 10^{-4}, \dots, 40 \times 10^{-4} \}$,
used to define $\ntime(\nu)$ such that $t_{\max}^{\obs} \approx \nu \,
t_{\max}^{\pred}$, we defined~$40$ instances of
problem~(\ref{lslost}). (The number of observations is~$\nobs \approx
0.2 \ntime(\nu) (\ndiv+1)$; $\nu=10^{-4}$ corresponds to $\nobs=427$,
while $\nu= 4 \times 10^{-3}$ corresponds to $\nobs=14{,}460$.) For
each instance, a solution satisfying~(\ref{stopcrit}) was computed
considering~$36$ different tolerances $\epsilon \in \{ 7.5 \times
10^{-13}, 5 \times 10^{-13}, 2.5 \times 10^{-13}, \dots, 10^{-4}
\}$. For each combination $(\nu,\epsilon)$, we obtained a solution
$\bar \xi(\nu,\epsilon)$, that is said to be acceptable
if~(\ref{accep}) holds. Figure~\ref{fig1} displays, as a function
of~$\nu$ and~$\epsilon$, the value of the prediction error~$\eta(\bar
\xi(\nu,\epsilon))$ defined in~(\ref{accep}). In the figure, cold
colors (blue, cyan, and green) correspond to solutions that are not
acceptable; while hot colors (yellow, orange, red, and dark red)
correspond to acceptable solutions. The figure shows (on the left)
that acceptable solutions were not found when the number of
observations~$\nobs(\nu)$ was smaller than the number of unknowns
$\ndiv=500$. When the number of observations is larger than the number
of unknowns, acceptable solutions are only found when $\epsilon \leq
10^{-9}$.

\subsubsection{Choice of the subproblems' dimension} \label{nr2}

We now consider an instance of problem~(\ref{lslost}) with
$\ndiv=500$, $x_{\min}=0$, $\Delta x=6$, $\ntime=10$, $t_{\min}=0$,
$\Delta t = 0.1$, and $S$ randomly chossen in such a way that $\nobs =
|S| \approx 0.1 (2 \ntime (\ndiv+1))$, i.e.\ assuming that
approximately 90\% of the observations are not available. The choice
$\ntime=10$ combined with $\Delta t=0.1$ means that observations are
collected at intervals of $0.1$ seconds during $1$ second; and since
we are assuming that 90\% of the observations will not be available,
this means that there will be $\nobs \approx 0.1 \times 2 \times 10
\times (\ndiv + 1) = 2 (\ndiv+1) > \ndiv$ observations available. We
aim to find solutions to this instance satisfying~(\ref{stopcrit})
with $\epsilon=10^{-9}$ that, for this instance, corresponds to
$\ftarget \approx 1.9633 \times 10^{-5}$. Due to analysis in the
previous paragraph, it is expected the computed solution to be
acceptable according to~(\ref{accep}); so the solution can be used to
make predictions for the next $3{,}559$ seconds.

\begin{figure}[ht!]
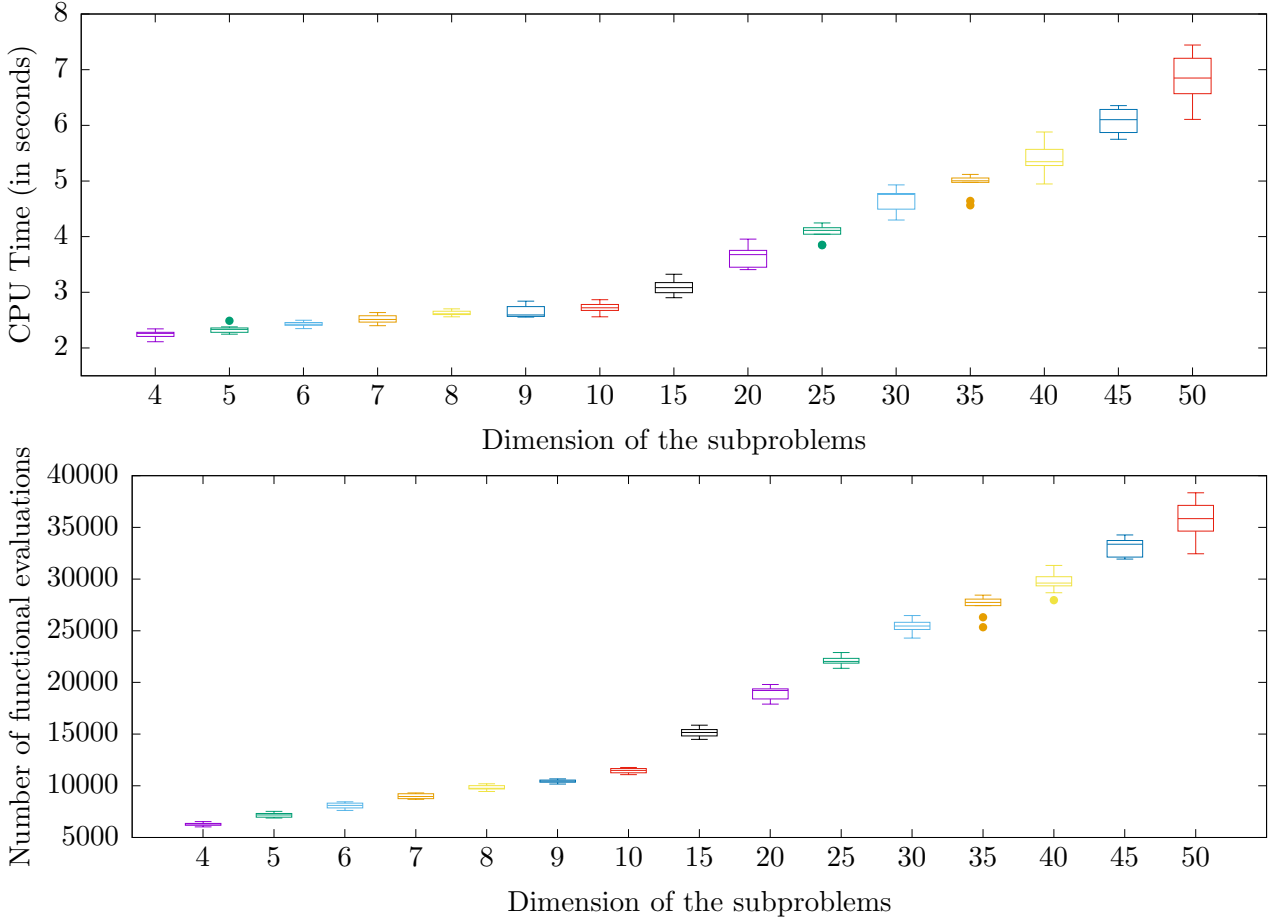

\begin{center}
\begin{tabular}{c}
\input{bmsesemfig2a.tex} \\
\input{bmsesemfig2b.tex}
\end{tabular}
\end{center}
\caption{Boxplots of performance metrics of Algorithm~\ref{algo}.1
  with the affine-subspaces-based reduction strategy applied to the
  instance with $\ndiv=500$ varying the subproblems' dimension
  $n_{\mathrm{red}}$.}
\label{fig2}
\end{figure}

\begin{figure}[ht!]
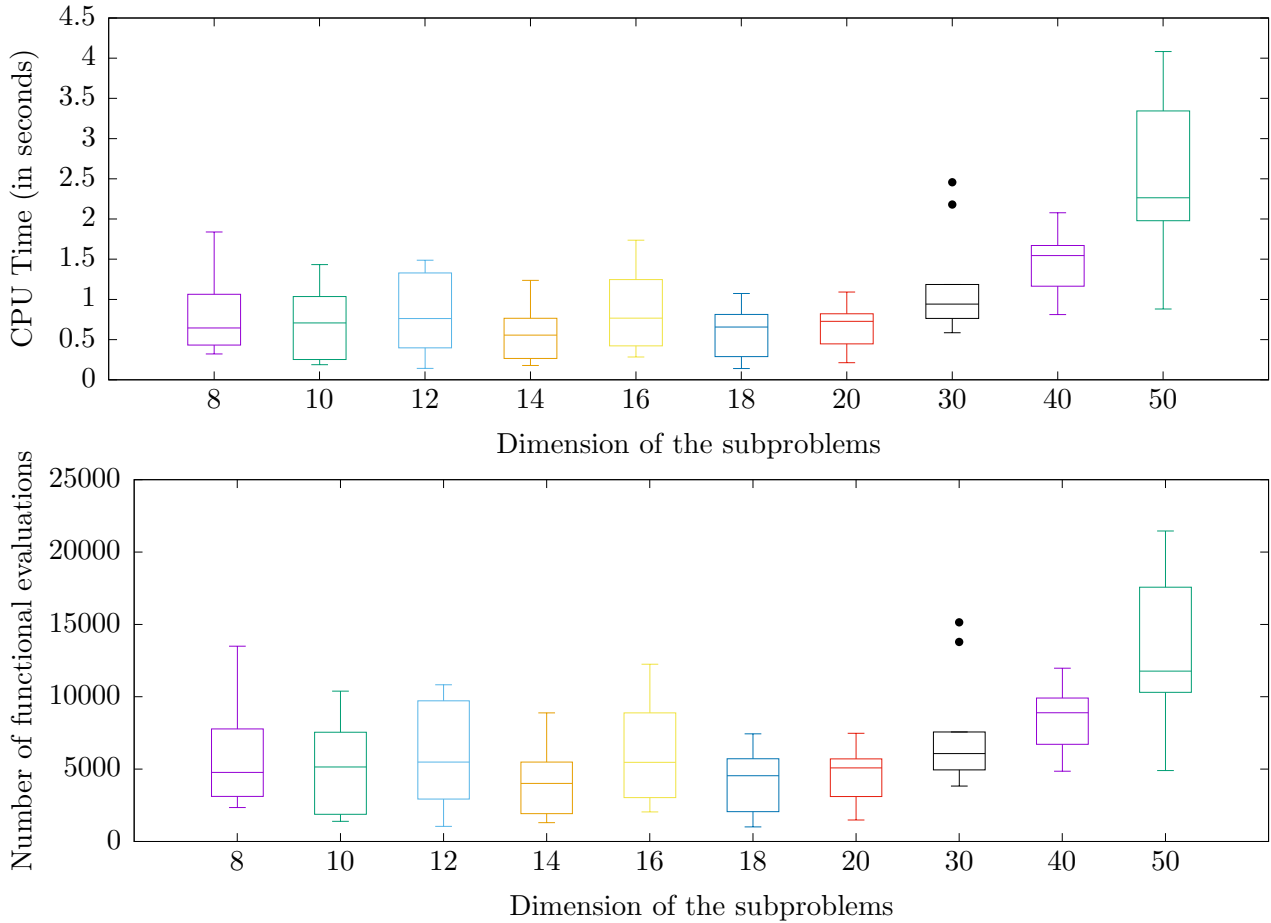

\begin{center}
\begin{tabular}{c}
\input{bmsesemfig3a.tex} \\
\input{bmsesemfig3b.tex}
\end{tabular}
\end{center}
\caption{Boxplots of performance metrics of Algorithm~\ref{algo}.1
  with the linear-interpolation-based reduction strategy applied to
  the instance with $\ndiv=500$ varying the subproblems' dimension
  $n_{\mathrm{red}}$.}
\label{fig3}
\end{figure}

The instance in the previous paragraph will be used to observe the
behavior of two variants of Algorithm~\ref{algo}.1, with
affine-subspaces-based and with linear-interpolation-based
subproblems, under variations of the subproblems' dimension
$\nred$. Each variation of Algorithm~\ref{algo}.1 uses, at every
iteration, the same reduction strategy and the same subproblem's
dimension. As mentioned in the previous paragraph, $\ftarget \approx
1.9633 \times 10^{-5}$; while the initial guess is always
$x^0=0$. Figures~\ref{fig2} and~\ref{fig3} show the results. Since
both reduction strategies have a random component, the instance was
solved ten times for each considered value
of~$\nred$. Figure~\ref{fig2} shows boxplots of two performance
measures (CPU time and number of functional evaluations) of
Algorithm~\ref{algo}.1 with affine-subspace-based subproblems and
$\nred \in \{ 4, 5, \dots, 9 \} \cup \{10, 15, \dots, 50\}$. The
boxplots show that the efficiency of the method is inversely
proportional to the size of the subproblems. Thus, it is worth
noticing that with $\nred \in \{ 2, 3\}$ the performance measures
present a large standard deviation and some outliers, while the method
fails a few times, characterizing a situation in which the method has
difficulties in improving the current approximation to a solution by
inspecting a very small search space. Figure~\ref{fig3} shows boxplots
of two performance measures (CPU time and number of functional
evaluations) of Algorithm~\ref{algo}.1 with linear-interpolation-based
subproblems and $\nred \in \{ 8, 10 ,12, 14, 16, 18, 20, 30, 40, 50
\}$. The boxplots show a uniform performance of the method for $\nred
\leq 20$; while, for $\nred >20$, the efficiency decreases when
$\nred$ increases.

\subsubsection{Influence of the acceleration scheme} \label{nr3}

\begin{figure}[ht!]
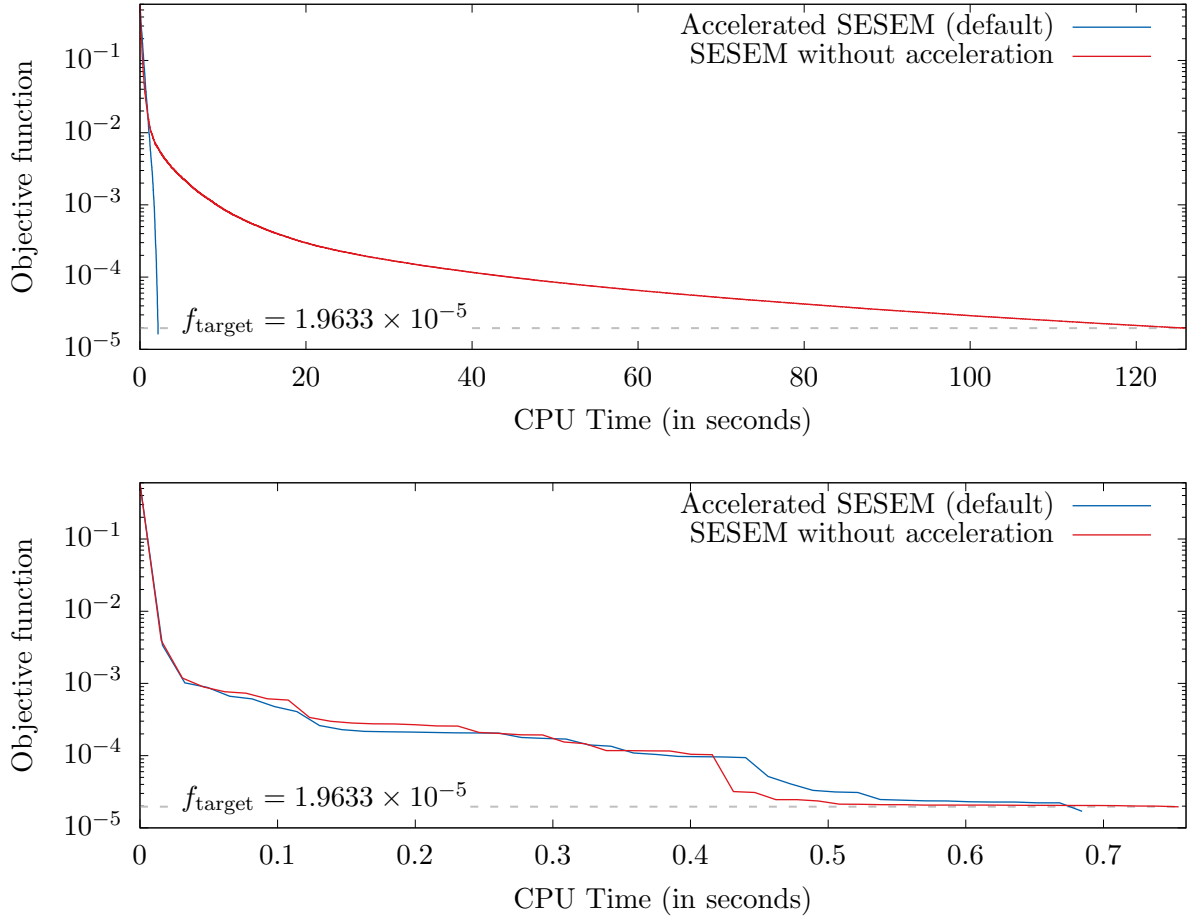

\begin{center}
\begin{tabular}{c}
\input{bmsesemfig4a.tex}\\
\input{bmsesemfig4b.tex}
\end{tabular}
\end{center}
\caption{Influence of the acceleration in the performance of
  Algorithm~\ref{algo}.1 with affine-subspaces-based subproblems (top)
  and with linear-interpolation-based subproblems (bottom) when
  applied to the instance with $\ndiv=500$.}
\label{fig4}
\end{figure}

Still considering the same instance, we now analyze the influence of
the acceleration in the performance of Algorithm~\ref{algo}.1 with
affine-subspaces-based subproblems ($\nred=4$) and with
linear-interpolation-based subproblems ($\nred=20$). Figure~\ref{fig4}
shows the results. The figure shows that, when the affine-subspaces
reduction strategy is considered, the acceleration improves the
efficiency of the method in approximately two orders of magnitude;
while it appears to have no relevant effect in combination with the
linear-interpolation-based reduction strategy; although it appears to
speed up the convergence of the method in its final iterations.

\subsubsection{Solving larger instances} \label{nr4}

\begin{figure}[ht!]
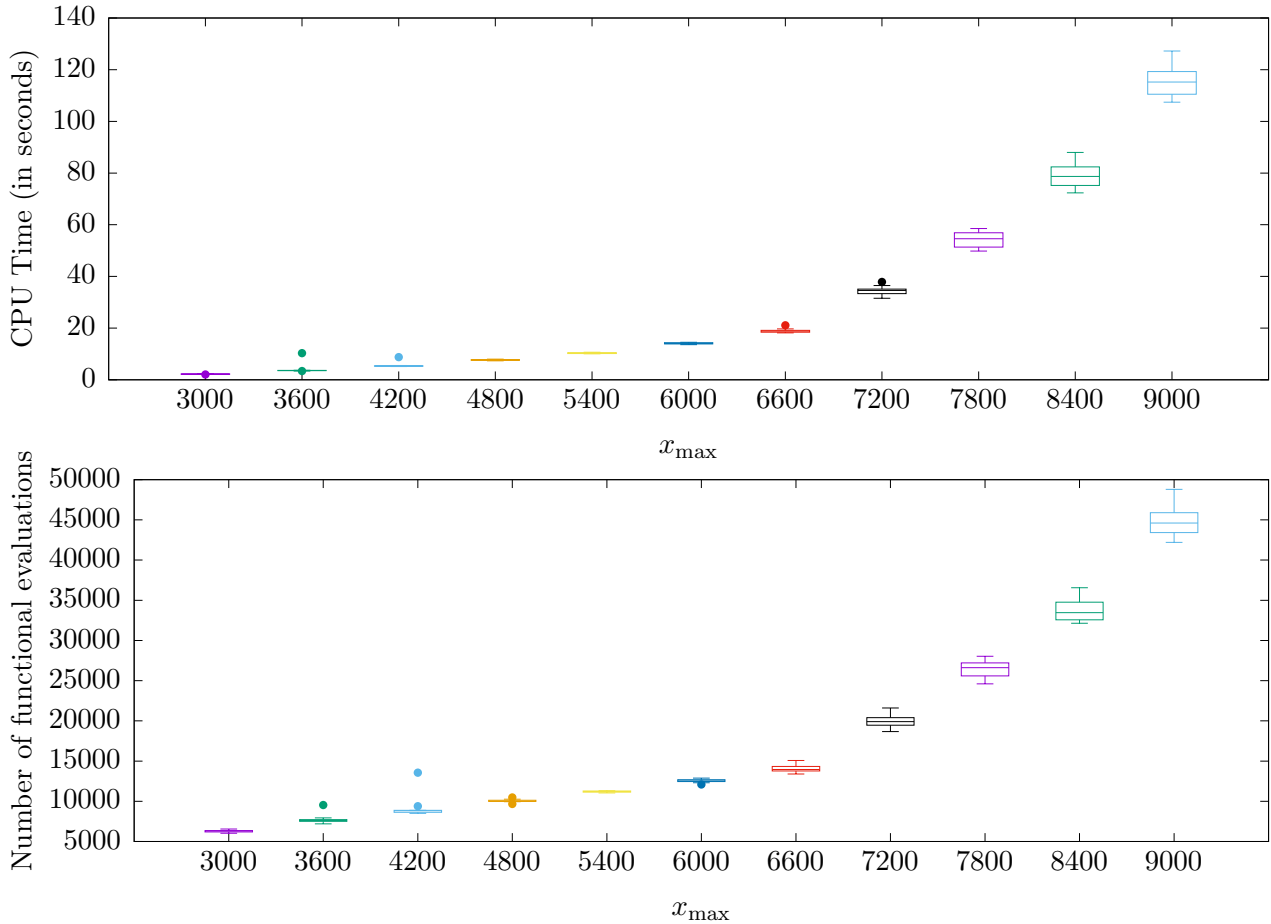

\begin{center}
\begin{tabular}{c}
\input{bmsesemfig5a.tex} \\
\input{bmsesemfig5b.tex}
\end{tabular}
\end{center}
\caption{Boxplots of performance metrics of Algorithm~\ref{algo}.1
  with affine-subspaces-based subproblems ($n_{\mathrm{red}}=4$)
  applied to instances of increasing size with $x_{\max} \in \{
  3{,}000, 3{,}600, \dots, 9{,}000 \}$.}
\label{fig5}
\end{figure}

\begin{figure}[ht!]
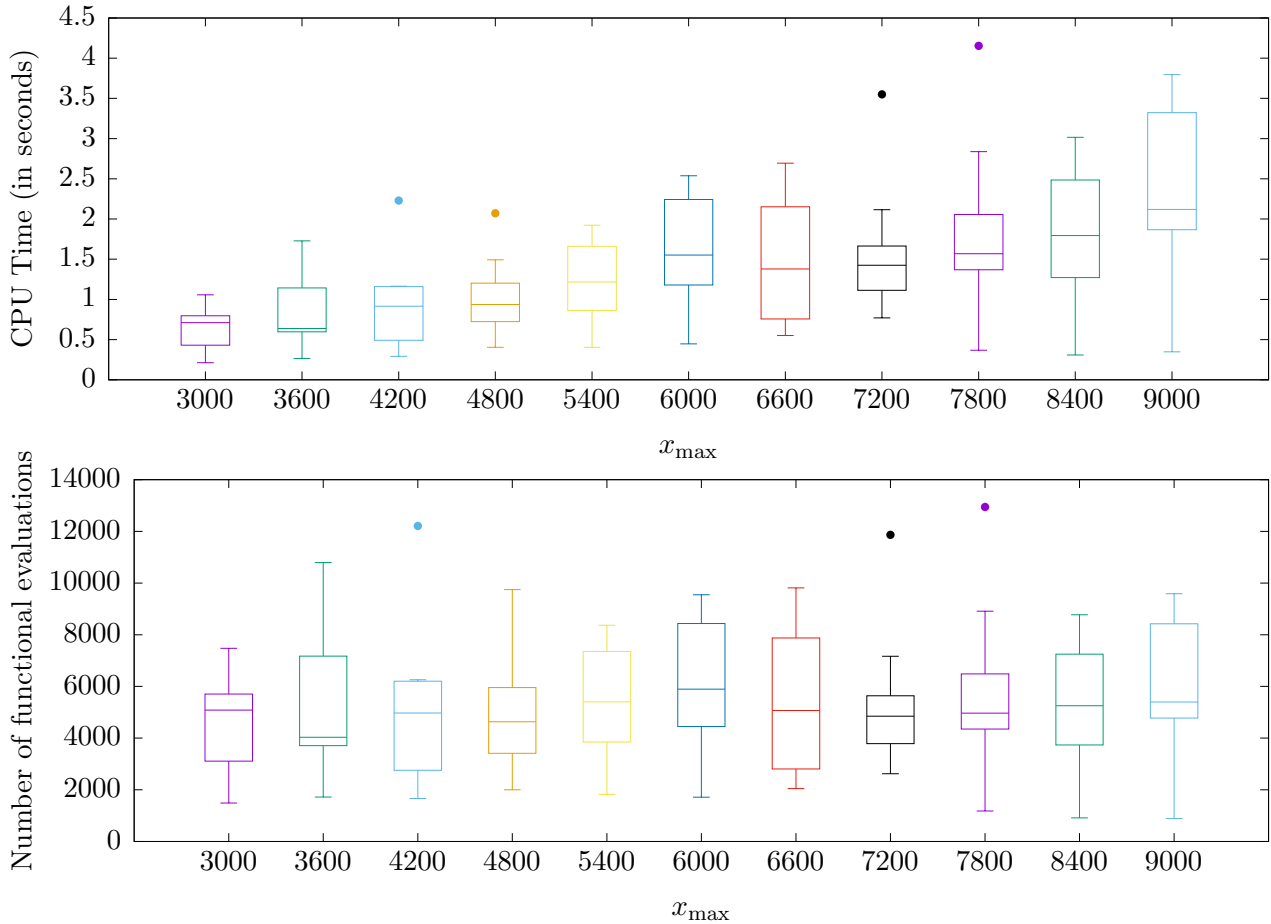

\begin{center}
\begin{tabular}{c}
\input{bmsesemfig6a.tex} \\
\input{bmsesemfig6b.tex}
\end{tabular}
\end{center}
\caption{Boxplots of performance metrics of Algorithm~\ref{algo}.1
  with linear-interpolation-based subproblems ($n_{\mathrm{red}}=20$)
  applied to instances of increasing size with $x_{\max} \in \{
  3{,}000, 3{,}600, \dots, 9{,}000 \}$.}
\label{fig6}
\end{figure}

We now consider a set of instances exactly as the one already
described but with $\ndiv \in \{ 500, 600, \dots$, $1{,}500
\}$. (These values correspond to $x_{\max} = 3{,}000, 3{,}600,
4{,}200, \dots, 9{,}000$, respectively.) Table~\ref{tab1} presents the
performance of Algorithm~\ref{algo}.1 with affine-subspaces-based
subproblems ($\nred=4$) and with linear-interpolation-based
subproblems ($\nred=20$). As before, the initial guess~$x^0$ is always
the origin, $\ftarget$ corresponds to the value of the right-hand-side
in~(\ref{stopcrit}) with $\epsilon=10^{-9}$. In the
table, $\|F(\bar \xi)\|_2^2$ corresponds to the left-hand-side
in~(\ref{stopcrit}), i.e.,
\[
\|F(\bar \xi)\|^2 = \sum_{\{ (i,j,k) \in S \}}
\left( y(\bar \xi,t_i,x_j,k) - y_{ijk}^{\obs} \right)^2,
\]
\#it stands for the number of iterations, \#fcnt stands for the number
of functional evaluations, and Time stands for the CPU time in
seconds. Since the method is run ten times per instance, values in the
table correspond to averages. In addition, for the CPU time, the
standard deviation is also presented in the table; and boxplots are
given in Figures~\ref{fig5} and~\ref{fig6}. A comparison between
Figures~\ref{fig5} and~\ref{fig6} makes it clear that the cost of the
affine-subspaces strategy grows together with the size of the
instances; while the linear-interpolation strategy appears to absorve
the cost of increasing sizes by incorporation some knowledge of the
problem's solution.

\begin{table}[ht!]
\begin{center}
\begin{tabular}{|rr|crrrc|crrrc|}
\hline
\multicolumn{1}{|c}{\multirow{3}{*}{$\ndiv$}} &
\multicolumn{1}{c|}{\multirow{3}{*}{$\nobs$}} &
\multicolumn{5}{c|}{SESEM with affine subspaces ($n_{\mathrm{red}}=4$)} &
\multicolumn{5}{c|}{SESEM with linear splines ($n_{\mathrm{red}}=20$)}\\
\cline{3-12}
&& \multirow{2}{*}{$\|F(\bar \xi)\|^2$} & \multicolumn{1}{c}{\multirow{2}{*}{\#it}} &
\multicolumn{1}{c}{\multirow{2}{*}{\#fcnt}} & \multicolumn{2}{c|}{Time} &
\multirow{2}{*}{$\|F(\bar \xi)\|^2$} & \multicolumn{1}{c}{\multirow{2}{*}{\#it}} &
\multicolumn{1}{c}{\multirow{2}{*}{\#fcnt}} & \multicolumn{2}{c|}{Time} \\
&&&&& avg & stdev &&&& avg & stdev\\
\hline
\hline
    500 & 1{,}058 & 1.86e-05 &     464 &  6{,}293 &   2.19 & 0.07 & 1.87e-05 & 38 & 4{,}598 & 0.64 & 0.29\\
    600 & 1{,}257 & 2.10e-05 &     568 &  7{,}765 &   4.25 & 2.14 & 1.99e-05 & 39 & 5{,}094 & 0.81 & 0.46\\
    700 & 1{,}471 & 2.63e-05 &     669 &  9{,}265 &   5.65 & 1.11 & 2.49e-05 & 42 & 5{,}109 & 0.93 & 0.54\\
    800 & 1{,}669 & 3.06e-05 &     740 & 10{,}051 &   7.65 & 0.21 & 2.89e-05 & 43 & 4{,}984 & 1.03 & 0.48\\
    900 & 1{,}908 & 3.53e-05 &     829 & 11{,}192 &  10.32 & 0.23 & 3.17e-05 & 45 & 5{,}514 & 1.25 & 0.51\\
1{,}000 & 2{,}106 & 3.83e-05 &     921 & 12{,}538 &  14.07 & 0.30 & 3.72e-05 & 45 & 6{,}021 & 1.59 & 0.65\\
1{,}100 & 2{,}312 & 4.26e-05 & 1{,}024 & 14{,}089 &  19.01 & 0.86 & 4.00e-05 & 41 & 5{,}304 & 1.45 & 0.75\\
1{,}200 & 2{,}484 & 4.66e-05 & 1{,}360 & 19{,}924 &  34.46 & 1.87 & 4.46e-05 & 46 & 5{,}364 & 1.58 & 0.80\\
1{,}300 & 2{,}677 & 5.05e-05 & 1{,}713 & 26{,}404 &  54.18 & 3.29 & 4.88e-05 & 46 & 5{,}802 & 1.84 & 1.03\\
1{,}400 & 2{,}885 & 5.45e-05 & 2{,}124 & 33{,}796 &  78.93 & 5.04 & 4.82e-05 & 43 & 5{,}227 & 1.79 & 0.91\\
1{,}500 & 3{,}090 & 5.83e-05 & 2{,}725 & 44{,}897 & 115.64 & 6.70 & 5.05e-05 & 48 & 5{,}916 & 2.34 & 1.01\\
\hline
\end{tabular}
\end{center}
\caption{Performance of Algorithm~\ref{algo}.1 with
  affine-subspaces-based subproblems ($\nred=4$) and with
  linear-interpolation-based subproblems ($\nred=20$) applied to
  instances of increasing size with $x_{\max} \in \{ 3{,}000, 3{,}600,
  \dots, 9{,}000 \}$.}
\label{tab1}
\end{table}

\subsubsection{Comparison with BOBYQA and DFBOLS} \label{nr5}

This section ends presenting the performance of BOBYQA and
DFBOLS\footnote{Provided by Hongchao Zhang on January 11th, 2021.}
applied to the same instances of Table~\ref{tab1}. Aiming a fair
comparison, both methods were modified to stop as soon as they reach a
solution~$\xi$ satisfying $\| F(\xi)\|_2^2 \leq
\ftarget$. Table~\ref{tab2} shows the results. Figures in the table
show that both variants of Algorithm~\ref{algo}.1 outperforms BOBYQA
and DFBOLS by several orders of magnitude when the CPU time is
considered as performance measure. While DFBOLS is the most time
consuming method, it is the most efficient if the number of functional
evaluations is considered. It is worth noticing that the comparison
between the behaviors of the considered methods is restricted to their
application to the problem under consideration.

\begin{table}[ht!]
\begin{center}
\begin{tabular}{|rr|crr|crr|}
\hline
\multicolumn{1}{|c}{\multirow{2}{*}{$\ndiv$}} &
\multicolumn{1}{c|}{\multirow{2}{*}{$\nobs$}} &
\multicolumn{3}{c|}{BOBYQA} &
\multicolumn{3}{c|}{DFBOLS}\\
\cline{3-8}
&& $\|F(\bar \xi)\|^2$ & \multicolumn{1}{c}{\#fcnt} & \multicolumn{1}{c|}{Time} &
   $\|F(\bar \xi)\|^2$ & \multicolumn{1}{c}{\#fcnt} & \multicolumn{1}{c|}{Time} \\
\hline
\hline
    500 & 1{,}058 & 1.96e-05 &  7{,}593 &      278.68 & 2.63e-07 & 1{,}006 &  1{,}907.91 \\
    600 & 1{,}257 & 2.34e-05 &  8{,}176 &      511.07 & 2.39e-07 & 1{,}206 &  5{,}679.65 \\
    700 & 1{,}471 & 2.73e-05 & 10{,}601 &  1{,}057.75 & 5.16e-07 & 1{,}406 & 13{,}413.81 \\
    800 & 1{,}669 & 3.14e-05 & 15{,}353 &  1{,}948.76 & 1.16e-07 & 1{,}606 & 31{,}103.48 \\
    900 & 1{,}908 & 3.58e-05 & 12{,}901 &  2{,}133.89 & --& --& $>$10h \\
1{,}000 & 2{,}106 & 3.92e-05 & 21{,}674 &  4{,}376.88 & --& --& --\\
1{,}100 & 2{,}312 & 4.32e-05 & 18{,}916 &  5{,}967.14 & --& --& --\\
1{,}200 & 2{,}484 & 4.67e-05 & 22{,}419 & 11{,}584.94 & --& --& --\\
1{,}300 & 2{,}677 & 5.06e-05 & 26{,}074 & 16{,}742.44 & --& --& --\\
1{,}400 & 2{,}885 & 5.46e-05 & 34{,}286 & 26{,}692.62 & --& --& --\\
1{,}500 & 3{,}090 & 5.84e-05 & 29{,}473 & 26{,}325.37 & --& --& --\\
\hline
\end{tabular}
\end{center}
\caption{Performance of BOBYQA and DFBOLS applied to instances of
  increasing size with $x_{\max} \in \{ 3{,}000, 3{,}600, \dots,
  9{,}000 \}$.}
\label{tab2}
\end{table}

\section{Final remarks} \label{conclusions}

In this paper, we presented a general scheme under which globally
convergent derivative-free algorithms for nonlinear least squares with
sequential secant acceleration can be defined. Our main motivation was
the estimation of parameters in hydraulic models governed by partial
differential equations. The non-availability of derivatives come from
the fact that these models may be computed by ``partially black-box''
codes and the possible uncertainty of function evaluations motivated
by the lack of precise geometrical parameters during the estimation
process.

Algorithms based on interpolating quadratic models like
BOBYQA~\cite{bobyqa} (see, also, \cite{csv}) use to be effective for
this type of optimization problems. However, big costs associated with
model building and its minimization make it necessary to employ
schemes in which the number of variables is not very large. Partial
minimization over random affine subspaces is an adequate
dimension-reduction
procedure~\cite{cartisroberts2021,yaxiang1,yaxiang2}; and we showed
that acceleration based on the sequential secant framework is
effective to increase the performance of that approach. In addition,
we developed a new reduction procedure based on variable linear
interpolation in which the variables of subproblems are a set of
independent variables with coordinates corresponding to (also)
variable nodes. The effectiveness of this new approach is associated
with the structure of the variables of the problem. If, in the
original underlying problem, the unknown is a continuous function that
depends on a single variable, the one-dimensional interpolatory scheme
tends to be quite effective. This is the case of our problem of
estimating the Manning coefficients, which, as a consequence, does not
need acceleration to obtain the best possible results. In more
complicated cases, the ``true'' unknown of the problem may be a
continuous function of~2, 3, or more variables. In this case, our
variable-node interpolation scheme should be conveniently adapted by
means of incorporation of multi-dimensional interpolation devices. In
the present work, the new dimensionality-reduction scheme, as well as
the acceleration process, were apllied in connection with the
derivative-free general-purpose solver BOBYQA~\cite{bobyqa}. In
exactly the same way, both features can be used in connection with
derivative-free least-squares methods such as the ones introduced
in~\cite{cartisroberts2019,cartisroberts2021,zcs}.

\end{document}